\magnification = 1200 \baselineskip = 15 pt

\input psfig.sty

\def \P {{\bf P}}
\def \E {{\bf E}}
\def \Vt {\tilde{V}}
\def \V {V^f}

\def \p {\psi^{nb}}

\def \pv {\pi^{V}}

\def \Ma {M^{*}_{\tau_d}}

\def\sqr#1#2{{\vcenter{\vbox{\hrule height.#2pt\hbox{\vrule width.#2pt height#1pt \kern#1pt\vrule width.#2pt}\hrule height.#2pt}}}}

\def \square{\hfill\mathchoice\sqr56\sqr56\sqr{4.1}5\sqr{3.5}5}

\def \sect#1{\bigskip \noindent {\bf  #1} \medskip}
\def \subsect#1{\medskip \noindent{\it #1} \medskip}
\def \th#1#2{\medskip \noindent {\bf  Theorem #1.}   \it #2 \rm}

\def \cor#1#2{\medskip \noindent {\bf  Corollary #1.}   \it #2 \rm}
\def \pf {\noindent  {\it Proof}.\quad }

\def \ex#1{\medskip \noindent {\bf  Example #1.}}

\centerline{\bf  Minimizing the Lifetime Shortfall or Shortfall at
Death} \bigskip \bigskip

\centerline{Erhan Bayraktar} \medskip

\centerline{Virginia R. Young} \bigskip \bigskip

\centerline{Department of Mathematics, University of
Michigan}\medskip \centerline{Ann Arbor, Michigan, 48109} \bigskip
\bigskip


\vfill \eject

\centerline{\bf  Minimizing the Lifetime Shortfall or Shortfall at
Death} \bigskip

\noindent{\bf  Abstract:}  We find the optimal investment strategy
for an individual who seeks to minimize one of four objectives:
(1) the probability that his wealth reaches a specified ruin level
{\it before} death, (2) the probability that his wealth reaches
that level {\it at} death, (3) the expectation of how low his
wealth drops below a specified level {\it before} death, and (4)
the expectation of how low his wealth drops below a specified
level {\it at} death.  Young (2004) showed that under criterion
(1), the optimal investment strategy is a heavily leveraged
position in the risky asset for low wealth.

In this paper, we introduce the other three criteria in order to
reduce the leveraging observed by Young (2004).  We discovered
that surprisingly the optimal investment strategy for criterion
(3) is {\it identical} to the one for (1) and that the strategies
for (2) and (4) are {\it more} leveraged than the one for (1) at
low wealth.  Because these criteria do not reduce leveraging, we
completely remove it by considering problems (1) and (3) under the
restriction that the individual cannot borrow to invest in the
risky asset.

\medskip


\medskip


\sect{1. Introduction}

We study an individual investment problem by using optimal
stochastic control.  The study of the investment problems faced by
individuals is justified because a significant financial crisis is
looming: it is projected that retired Americans' living expenses
will exceed their financial resources by \$400 billion over the
ten-year period 2020-2030 (VenDerhei and Copeland, 2003). This
shortfall is due to the increased longevity of our aging
population, changes in Social Security, inadequate private
retirement savings, and the continuing switch from defined benefit
plans to defined contribution plans, such as 401(k)s, which
transfers the investment and longevity risk from the employer to
the individual.

We consider the problem of how an individual should invest her
wealth in a risky financial market in order to minimize her (1)
probability of lifetime ruin, that is, the probability of running
out of money before dying, (2) probability of ruin at death, that
is, the probability of running out of money at the time of death
(this objective can be used by people with bequest motives), (3)
expected lifetime shortfall, and (4) expected shortfall at death.
The latter two criteria are the counterparts of (1) and (2) in
which the individual is penalized by the amount of loss. In (1)
and (2) the penalty she gets is constant regardless how low the
wealth becomes with respect to the ruin threshold, the value of
wealth at which the individual considers herself ruined.

As employers shift from defined benefit plans to defined
contribution plans, the problem of outliving one's wealth becomes
more and more relevant to retirees because the income coming from
guaranteed sources is projected to drop significantly; see, for
example, Parikh (2003). We determine the optimal investment
strategy of an individual who targets a given rate of consumption
under each criterion and show how the strategies compare to each
other.  We assume that the rate of consumption (either nominal or
real) is net of any income the retiree receives from pension
plans, such as Social Security or a defined benefit plan. In
finding the optimal strategy, we take into account that the time
of death is random. This assumption differs from the one usually
assumed by financial planners and common retirement planning
software in that they generally assume a specific age of death.

We focus on minimizing the expected maximum lifetime shortfall and
the shortfall at death. One might argue that the former penalty
criterion is too severe because if an individual were to have a
negative wealth at some point in life and later get out of debt,
then the person should be as ``happy" as someone who was never in
debt but possesses the same current wealth. This argues that the
latter objective function is more appropriate than the former;
that is, all that matters to the individual is whether she is in
debt at the time of death. On the other hand, if someone's wealth
were to become negative, then that could affect that person's
ability to borrow money in the future. Therefore, under this view,
the first objective of minimizing the expected maximum shortfall
during life is more appropriate.

The most common optimization criterion used in the mathematical
finance literature is to maximize one's expected discounted
utility of consumption and bequest; see, for example, Merton
(1992) and Karatzas and Shreve (1998, Chapter 3).  Also, see
Zari-phopoulou (1999, 2001) for helpful summaries of the work to
date in this area. The goal of maximizing expected discounted
utility of consumption and bequest may be difficult to implement
because it depends on a subjective utility function for
consumption and bequest.  Minimizing the criteria we suggest might
prove easier for individuals to understand because these criteria
are arguably more objective. However, one should note there is a
{\it correspondence} between the utility maximization problem and
the ruin minimization problems when the utility function in
question is HARA. But HARA is the only utility function when one
finds a correspondence (see Bayraktar and Young (2007b) for
details).

Recently, a variety of papers in the risk and portfolio management
literature revitalized the Roy (1952) ``Safety-First rule" and
applied the concept in the context of maximizing the probability
of achieving certain investment goals before ruin.  For example,
Browne (1995, 1997, 1999a,b,c) derived the optimal stretegy for a
portfolio manager who is interested in maximizing the probability
of reaching a safe level before ruin or minimizing the time
expected time it takes to reach a goal, the probability of beating
a stochastic benchmark, and the probability of reaching a goal by
a deadline.  Also, researchers have begun to study the problem of
optimal investment to minimize the probability that an individual
runs out of money before dying.  See, for example, Milevsky, Ho,
and Robinson (1997), Milevsky and Robinson (2000), and Young
(2004).  In insurance mathematics, the criteria of minimizing the
probability of ruin was used to find the optimal reinsurance; see,
for example, Schmidli 2001, Taksar and Markussen (2003), and
Promislow and Young (2005).

Young (2004) studied the problem of finding the investment
strategy to minimize the probability of lifetime ruin. She
discovered that the individual leveraged her wealth when wealth
approaches zero; that is, when wealth was low, the individual
borrowed money to invest in the risky asset. In order to reduce or
eliminate this leveraging, Bayraktar and Young (2007a) imposed
borrowing constraints on the investment strategy. They first
considered the case for which the individual is not allowed to
borrow money in order to invest in the risky asset, and this
certainly eliminated the leveraging that Young (2004) observed.
They next considered the case for which the individual was allowed
to borrow money but at a rate higher than the one earned by the
riskless asset.  They discovered that the leveraging in this case
could be even worse than that in the case when the individual is
allowed to borrow at the rate earned by the riskless asset.

In this paper, we show that the leveraging effect observed in the
probability of lifetime ruin problem is not reduced by considering
alternative penalty functions.  We obtain two rather surprising
from our models.  (1)  We learn that the leveraging is {\it
exacerbated} by considering the probability of ruin at death and
the shortfall at death.  (2)  The optimal investment strategy for
the lifetime shortfall is {\it identical} to the optimal
investment strategy for the probability of lifetime ruin. This
result is surprising because we believed that by penalizing the
individual for the magnitude of her bankruptcy -- not for just
whether or not bankruptcy occurred -- we would be able to temper
her leveraging.

Our contributions to the applied probability literature are as
follows:  (1) By applying the Fenchel-Legendre transform to our
value functions, we are able to linearize the corresponding
non-linear HJB equations.  The resulting linear differential
equation is one with a free boundary, and we solve this
free-boundary problem for the case of minimizing the probability
of ruin at death.  The corresponding minimum ruin probability is
not explicitly available, but via we characterize it as the
inverse Fenchel-Legendre transform of the solution of the
free-boundary problem.  (2) In the problem of minimizing the
expectation of any function of lifetime minimum wealth, we also
rely on the Fenchel-Legendre transform to show that corresponding
optimal investment strategy is identical to the one for minimizing
the probability of lifetime ruin.  Therefore, the
optimally-controlled wealth processes for both problems are
identical.  From this observation, we develop a representation of
lifetime shortfall in terms of the probability of lifetime ruin;
see equation (3.16) below.     (3)  We provide verification
theorems that show that the value function is the unique solution
of an associated Hamilton-Jacobi-Bellman (HJB) equation.  Our
verification theorems are for a diffusion with killing, as well as
include an additional state variable for the minimum wealth
process.

The rest of the paper is organized as follows: In Section 2.1, we
introduce the financial market and summarize the results of Young
(2004) concerning the optimal investment strategy to minimize the
probability of lifetime ruin. In Section 2.2, we solve the problem
of minimizing probability of ruin at death.  We discover that the
optimal investment strategy for this criterion is always greater
than the one for the probability of lifetime ruin from Section
2.1. In Section 2.3, we remove leveraging entirely by prohibiting
borrowing.

In Section 3.1, we consider the problem of minimizing lifetime
expected shortfall. We provide a verification theorem for the
problem of minimizing a bounded, twice continuously differentiable
function of minimum wealth. We construct a solution to this
auxiliary problem. The value function of the lifetime shortfall is
shown to be the limit of a sequence of appropriately-defined
auxiliary problems. The optimal investment strategy for the
lifetime shortfall is {\it identical} to the one for the
probability of ruin. In Section 3.2, we consider the problem of
minimizing shortfall at death and show that the optimal investment
strategy is greater than the one for the probability of lifetime
ruin. In Section 3.3, we eliminate leveraging in the lifetime
shortfall problem by prohibiting borrowing. Section 4 concludes
the paper.

\sect{2. Probability of Lifetime Ruin and Ruin at Death}

In Section 2.1, we present the financial market and briefly review
the problem of minimizing the probability of lifetime ruin.  We
provide an explicit expression for the optimal investment in the
risky asset and point out the leveraging at low wealth.  In
Section 2.2, we consider the problem of minimizing the probability
that wealth at death lies below a certain level.  We show that the
optimal investment in the risky asset under the latter problem is
greater than under the former; therefore, the leveraging has
increased, not decreased.  Therefore, because the change in value
function did not reduce the leveraging, in Section 2.3, we remove
it completely by prohibiting borrowing to invest in the risky
asset.

\subsect{2.1. Financial Market and Probability of Lifetime Ruin}

In this section, we first present the financial ingredients that
make up the individual's wealth throughout this paper, namely,
consumption, a riskless asset, and a risky asset. We, then,
determine the minimum probability of lifetime ruin. The individual
consumes at a constant continuous rate $c$. We assume that the
individual invests in a riskless asset whose price at time $t$,
$X_t$, follows the process $dX_t = rX_t dt, X_0 = x > 0$, for some
fixed rate of interest $r > 0$. Also, the individual invests in a
risky asset whose price at time $t$, $S_t$, follows geometric
Brownian motion given by

$$\left\{ \eqalign{dS_t &= \mu S_t dt + \sigma S_t dB_t, \cr S_0
&= S > 0,} \right. \eqno(2.1)$$ \noindent in which $\mu > r$,
$\sigma > 0$, and $B$ is a standard Brownian motion with respect
to a filtration of the probability space $(\Omega, {\cal F}, \P)$.
Let $W_t$ be the wealth at time $t$ of the individual, and let
$\pi_t$ be the amount that the decision maker invests in the risky
asset at that time.  It follows that the amount invested in the
riskless asset is $W_t - \pi_t$. Thus, wealth follows the process

$$\left\{ \eqalign{dW_t &= [rW_t + (\mu - r) \pi_t  - c(W_t)] dt +
\sigma \pi_t dB_t, \cr W_0 &= w.} \right. \eqno(2.2)$$

\noindent A process associated with this wealth process is the
minimum wealth process.  Let $M_t$ denote the minimum wealth of
the individual during $[0, t]$; that is,

$$M_t = \min\left[ \inf_{0 \le s \le t} W_s, M_0\right],
\eqno(2.3)$$

\noindent in which we include $M_0$ (possibly different from
$W_0$) to allow for the individual to have a financial past.

By ``outliving her wealth,'' or equivalently ``lifetime ruin,'' we
mean that the individual's wealth reaches some specific value $x <
c/r$ before she dies.  Note that $c/r$ is the ``safe'' level above
which the individual cannot ruin if she places at least $c/r$ in
the riskless asset.  Let $\tau_x$ denote the first time that
wealth equals $x < c/r$, and let $\tau_d$ denote the random time
of death of our individual.  We assume that $\tau_d$ is
exponentially distributed with parameter $\lambda$ (that is, with
expected time of death equal to $1/\lambda$); this parameter is
also known as the {\it hazard rate} of the individual.

Denote the minimum probability that the individual outlives her
wealth by $\psi(w, m; x)$, in which the argument $w$ indicates
that one conditions on the individual possessing wealth $w$ and
minimum wealth $m$ at the current time, and the parameter $x$
reminds us of the hitting level.  Thus, $\psi$ is the minimum
probability that $\tau_x < \tau_d$, in which one minimizes with
respect to admissible investment strategies $\pi$.  A strategy
$\pi$ is {\it admissible} if it is ${\cal F}_t$-progressively
measurable (in which ${\cal F}_t$ is the augmentation of
$\sigma(W_s: 0 \le s \le t)$) and if it satisfies the
integrability condition $\int_0^t \pi_s^2 \, ds < \infty$ almost
surely for all $t \ge 0$.

Note that the event that $\tau_d < \tau_x$ is the same event that
$M_{\tau_d} \le x$ with probability 1; thus, we can express $\psi$
as

$$\psi(w, m; x) = \inf_{\pi} \P \left[M_{\tau_d} \le x \vert W_0 =
w, M_0 = m \right]. \eqno(2.4)$$

\noindent Young (2004) showed that

$$\psi(w, m; x) = \cases{1, &if $m \le x$; \cr \left( {c - r w
\over c - r x} \right)^p, &if $x < m \le w < c/r$; \cr 0, &if $x <
m$ and $w \ge c/r$;}  \eqno(2.5)$$ \noindent with $$p = {1 \over
2r} \left[ (r + \lambda + \delta) + \sqrt{(r + \lambda + \delta)^2
- 4r \lambda} \right] > 1, \eqno(2.6)$$ and $$\delta = {1 \over 2}
\left( {\mu - r \over \sigma} \right)^2. \eqno(2.7)$$

Young (2004) also showed that the corresponding optimal investment
in the risky asset $\pi^\psi$ for $x < m \le w < c/r$ is given by

$$\pi^\psi(w) = {\mu - r \over \sigma^2} {1 \over p - 1} \left( {c
\over r} - w \right), \eqno(2.8)$$

\noindent a positive, decreasing, linear function of wealth,
independent of both $m$ and $x$.  Note that as wealth increases
towards $c/r$, the amount invested in the risky asset decreases to
zero. This makes sense because as the individual becomes
wealthier, she does not need to take on as much risk to achieve
her fixed consumption rate of $c$.  On the other hand, for wealth
small, the optimal amount invested in the risky asset is greater
than wealth; that is, the individual borrows money to invest in
the risky asset in order to avoid the greater risk of lifetime
ruin.

Browne (1997) considered a problem similar, but not equivalent, to
minimizing the probability of ruin.  He sought to maximize the
probability of hitting the safe level, $c/r$, before reaching zero
wealth; Pestien and Sudderth (1985) solved a related problem for a
general diffusion and showed that the optimal investment strategy
is determined by maximizing the ratio of the drift to the
variance.  (As an aside, their technique is not applicable if
$\lambda \ne 0$.)  However, Browne's problem does not have an
optimal investment strategy because the safe level is not
attainable.  Our problem of minimizing the probability of lifetime
ruin, as given in (2.4), is not equivalent to maximizing the
probability of reaching a certain level.  An optimal strategy
exists for (2.4), as given by (2.8), in contrast to Browne's case.

Browne (1997) obtained an $\epsilon$-optimal policy for his
problem that exhibits behavior similar to that given in (2.8),
that is, with leveraging at low wealth.  We believe that most
people with small wealth will not borrow to invest in a risky
asset to avoid ruin and that no credible financial advisor will
give such advice.  Therefore, in the next section, we modify the
objective function to try to reduce this leveraging at low wealth.

\subsect{2.2. Probability of Ruin at Death}

In this section, we modify the objective function in hopes that
the leveraging effect observed in (2.8) will be reduced.  One
might argue that the penalty that lifetime minimum wealth reaches
a certain level is too severe and that the individual is happy
enough as long as her wealth at death lies above that level,
regardless of the path that wealth follows between now and then.
However, we cannot simply minimize $\P[W_{\tau_d} \le x \vert W_0
= w]$ for all $w < c/r$ because we expect the minimum probability
to be convex, as is $\psi$ in (2.5) on the interval $(x, c/r)$.
Recall that there exists no convex function on $(-\infty, c/r)$
such that the function is bounded and decreasing.

Therefore, consider the following related problem:

$$\phi(w, m; x, M) = \inf_\pi \P \left[ \{W_{\tau_d} \le x\} \cup
\{ M_{\tau_d} \le M \} \vert W_0 = w, M_0 = m \right],
\eqno(2.9)$$

\noindent in which $M < x < c/r$.  The domain of $\phi$ with
respect to $w$ is effectively $[M, c/r)$ over which we expect
$\phi$ to be convex.  By setting $M$ to be a large negative
number, $\phi$ approximates what one might mean by the minimum
probability of ruin at death.

For $m \le M$, $\phi$ is identically 1, and for $m > M$ and $w \ge
c/r$, $\phi$ is identically 0.  For $M < m \le w < c/r$, $\phi$
solves the following Hamilton-Jacobi-Bellman (HJB) equation:

$$\left\{ \eqalign{& \lambda (\phi(w) - {\bf  1}_{\{w \le x\}}) =
(rw - c) \phi'(w) + \min_\pi \left[ (\mu - r) \pi \phi'(w) + {1
\over 2} \sigma^2 \pi^2 \phi''(w) \right], \cr & \phi(M) = 1,
\quad \phi(c/r) = 0,} \right. \eqno(2.10)$$

\noindent in which we drop the notational dependence of $\phi$ on
$m$ because $\phi$ is independent of $m$ if $m > M$.

The solution of (2.10) will be $C^1(M, c/r) \cap C^2((M, c/r) -
\{x\})$; therefore, because the solution is not smooth at $w = x$,
a corresponding verification theorem relies on the approximation
technique used in \O ksendal (1998, proof of Theorem 10.4.1).  Via
a verification theorem similar to the ones in Bayraktar and Young
(2007a) modified by this approximation technique, one can show
that the solution of the boundary-value problem in (2.10) is the
minimum probability in (2.9).    We proceed to solve (2.10).

\th{2.1} {For $m > M$ and $x \le w < c/r$, the function $\phi$ is
given by}

$$\phi(w) = \beta \left( 1 - {r \over c} w \right)^p,
\eqno(2.11)$$

\noindent{\it with $p$ as in $(2.6)$ and $\beta > 0$.  Thus, on
$[x, c/r)$, $\phi$ is a multiple of the probability of lifetime
ruin $\psi$.  It follows that on $[x, c/r)$, the optimal
investment strategy $\pi^\phi$ equals $\pi^\psi$ as given in
$(2.8)$.}

\medskip

\pf From (2.10), for $x \le w < c/r$, $\phi$ solves

$$\lambda \phi(w) = (rw - c) \phi'(w) + \min_{\pi} \left[ (\mu -
r) \pi \phi'(w) + {1 \over 2} \sigma^2 \pi^2 \phi''(w) \right],
\eqno(2.12)$$

\noindent with boundary condition $\phi(c/r) = 0$.  We also have
the boundary condition $\phi'(c/r) = 0$, which we demonstrate now.
Consider the minimum probability of lifetime ruin $\psi(w, m; x)$
given in (2.5).  Certainly, we have $0 \le \phi \le \psi$ because
the probability of ruining at level $M$ before dying or at level
$x > M$ at death is no greater than the probability of ruining at
level $x$ before dying.  Note that $\psi'(c/r) = 0$ from (2.5), in
which we differentiate with respect to $w$; therefore, $\phi'(c/r)
= 0$ because $\phi$ is wedged between 0 and $\psi$ as wealth
approaches $c/r$.

We hypothesize that $\phi$ is convex on $[x, c/r)$, and we
consider its Fenchel-Legendre transform $\tilde \phi$ defined by

$$\tilde \phi(y) = \min_w [\phi(w) + wy]. \eqno(2.13)$$

\noindent Note that we can recover $\phi$ from $\tilde \phi$ by

$$\phi(w) = \max_y [\tilde \phi(y) - wy]. \eqno(2.14)$$

\noindent The minimizing value of $w$ in (2.13) equals $I(-y) =
\tilde \phi'(y)$, in which $I$ is the inverse function of $\phi'$.
Therefore, the maximizing value of $y$ in (2.14) equals
$-\phi'(w)$.

Substitute $w = I(-y)$ in equation (2.12) to obtain

$$\lambda \tilde \phi(y) + (r - \lambda) y \tilde \phi'(y) -
\delta y^2 \tilde \phi''(y) = cy, \eqno(2.15)$$

\noindent in which $\delta$ is given in (2.7).  The general
solution of (2.15) is

$$\tilde \phi(y) = A_1 y^{B_1} + A_2 y^{B_2} + {c \over r} y,
\eqno(2.16)$$

\noindent in which $A_1$ and $A_2$ are constants to be determined,
and $B_1$ and $B_2$ are the positive and negative roots,
respectively, of

$$-\lambda - (r - \lambda + \delta) B + \delta B^2 = 0.
\eqno(2.17)$$

\noindent Thus,

$$B_1 = {1 \over 2 \delta} \left[ (r - \lambda + \delta) +
\sqrt{(r - \lambda + \delta)^2 + 4 \lambda \delta} \right] > 1,
\eqno(2.18)$$

\noindent and

$$B_2 = {1 \over 2 \delta} \left[ (r - \lambda + \delta) -
\sqrt{(r - \lambda + \delta)^2 + 4 \lambda \delta} \right] < 0.
\eqno(2.19)$$

\noindent Note that $B_1 = p/(p - 1)$.

Define $y_c = -\phi'(c/r) = 0$; that is, $\tilde \phi'(0) = c/r.$
From the definition of $\tilde h$ in (2.13) and from $\phi(c/r) =
0$, at $y = y_c = 0$, we have

$$\tilde \phi(0) = 0. \eqno(2.20)$$

\noindent It follows that $A_2 = 0$.  We can, then, recover $\phi$
from (2.16) and (2.14) and obtain the expression for $\phi$ in
(2.11).

Because $\phi$ is convex on $[x, c/r)$, the optimal policy
$\pi^\phi$ is given by the first-order necessary condition in
(2.12).  Thus, $\pi^\phi(w) = \pi^\psi(w)$ for $x \le w < c/r$.
$\square$

\medskip

Next, we consider the solution of (2.10) for $M < m \le w < x$.
On this domain, $\phi$ solves the system

$$\left\{ \eqalign{& \lambda (\phi(w) - 1) = (rw - c) \phi'(w) +
\min_\pi \left[ (\mu - r) \pi \phi'(w) + {1 \over 2} \sigma^2
\pi^2 \phi''(w) \right], \cr & \phi(M) = 1, \quad \phi(x)/\phi'(x)
= -(c/r - x)/p,} \right. \eqno(2.21)$$

\noindent in which we assume that $\phi$ and $\phi'$ are
continuous at $w = x$.  Again, we consider the Fenchel-Legendre
transform of $\phi$ and show how to solve for $\tilde \phi$
explicitly.  From (2.14), we can then determine $\phi$.

By substituting $w = I(-y)$ in (2.21), we obtain an equation
similar to (2.15) whose general solution equals

$$\tilde \phi(y) = D_1 y^{B_1} + D_2 y^{B_2} + {c \over r} y + 1,
\eqno(2.22)$$

\noindent in which $B_1$ and $B_2$ are given in (2.18) and (2.19),
respectively, and $D_1$ and $D_2$ are constants to be determined.
To determine $D_1$ and $D_2$, we rely on the boundary conditions
of $\phi$ at $w = M$ and $w = x$.

If we define $y_x = -\phi'(x) = \beta p (r/c) (1 - rx/c)^{p-1}$,
then we have

$$\tilde \phi(y_x) = {1 \over p} y_x \left[ {c \over r} + (p-1)x
\right], \eqno(2.23)$$

\noindent and

$$\tilde \phi'(y_x) = x.  \eqno(2.24)$$

\noindent Similarly, if we define $y_M = -\phi'(M)$, then we have

$$\tilde \phi(y_M) = 1 + M y_M, \eqno(2.25)$$

\noindent and

$$\tilde \phi'(y_M) = M. \eqno(2.26)$$

If we write these four equations in terms of the expression for
$\tilde \phi$ in (2.22), we obtain, respectively,

$$D_1 y_x^{B_1} + D_2 y_x^{B_2} + {c \over r} y_x + 1 = {1 \over
p} y_x \left[ {c \over r} + (p-1)x \right], \eqno(2.27)$$

$$D_1 B_1 y_x^{B_1 - 1} + D_2 B_2 y_x^{B_2 - 1} + {c \over r} = x,
\eqno(2.28)$$

$$D_1 y_M^{B_1} + D_2 y_M^{B_2} + {c \over r} y_M + 1 = 1 + M y_M,
\eqno(2.29)$$

\noindent and

$$D_1 B_1 y_M^{B_1 - 1} + D_2 B_2 y_M^{B_2 - 1} + {c \over r} = M.
\eqno(2.30)$$

Solve equations (2.29) and (2.30) for $D_1$ and $D_2$ to obtain

$$D_1 = {1 - B_2 \over B_1 - B_2} \left( M - {c \over r} \right)
y_M^{1 - B_1} < 0, \eqno(2.31)$$

\noindent and

$$D_2 = {B_1 - 1 \over B_1 - B_2} \left( M - {c \over r} \right)
y_M^{1 - B_2} < 0. \eqno(2.32)$$

\noindent Substitute these expressions for $D_1$ and $D_2$ into
(2.28) to get an equation for $y_x/y_M$:

$${B_1(1 - B_2) \over B_1 - B_2} \left( {c \over r} - M \right)
\left( {y_x \over y_M} \right)^{B_1 - 1} + {B_2(B_1 - 1) \over B_1
- B_2} \left( {c \over r} - M \right)  \left( {y_x \over y_M}
\right)^{B_2 - 1} = {c \over r} - x.  \eqno(2.33)$$

\noindent Equation (2.33) has a unique solution $y_x/y_M \in (0,
1)$ because (i) the left-hand side equals $c/r - M > c/r - x$ when
$y_x/y_M = 1$; (ii) as $y_x/y_M$ approaches 0, the left-hand side
goes to $-\infty$; and (iii) the left-hand side is increasing with
respect to $y_x/y_M$.

Once we have the solution to (2.33), then we can solve for $y_x$
as follows.  First, substitute the expressions for $D_1$ and $D_2$
into equation (2.27) to obtain

$$\eqalign{{1 \over y_x} = -{p - 1 \over p} \left( {c \over r} - x
\right) & + {1 - B_2 \over B_1 - B_2} \left( {c \over r} - M
\right) \left( {y_x \over y_M} \right)^{B_1 - 1} \cr & + {B_1 - 1
\over B_1 - B_2} \left( {c \over r} - M \right) \left( {y_x \over
y_M} \right)^{B_2 - 1}.} \eqno(2.34)$$

\noindent Substitute for $y_x/y_M$ in the right-hand side of
(2.34) then solve for $y_x$.  One technical point is that the
right-hand side is required to be positive; this is true, but we
omit the proof.  Now, $y_M = y_x/(y_x/y_M)$, and finally we get
$D_1$ and $D_2$ from (2.31) and (2.32), respectively.  Also, note
that $\beta$ in (2.11) is given by

$$\beta = {c \over rp} \left( 1 - {r \over c} x \right)^{1 - p}
y_x.  \eqno(2.35)$$

We summarize these results in the following theorem:

\th{2.2} {For $M < m \le w < x$, the function $\phi$ is given by
the inverse Fenchel-Legendre transform $(2.14)$ of $\tilde \phi$
in $(2.22),$ in which $D_1$ and $D_2$ are given by $(2.31)$ and
$(2.32),$ respectively.  $y_x/y_M$ is the unique solution of
$(2.33)$ in $(0, 1),$ and $y_x > 0$ is given by $(2.34)$.}
$\square$

\medskip

Next, for $M < m \le w < x$, we compare the optimal investment
strategy $\pi^\phi$ with $\pi^\psi$, as given in (2.8).  Assume
that the ruin level for the $\psi$-problem is less than or equal
to $M$.  Recall that for $m > M$ and $x \le w < c/r$, the two
strategies are identical.  For $M < m \le w < x$, $\pi^\phi(w) >
\pi^\psi(w)$ if and only if

$$-{\phi'(w) \over \phi''(w)} > {1 \over p - 1} \left( {c \over r}
- w \right).  \eqno(2.36)$$

\noindent By substituting $w = I(-y)$, in which $I$ is the inverse
of $\phi'$, we obtain that $\pi^\phi > \pi^\psi$ if and only if

$$-y \tilde \phi''(y) > {1 \over p - 1} \left( {c \over r} -
\tilde \phi'(y) \right), \quad y_x < y < y_M. \eqno(2.37)$$

\noindent By substituting into (2.37) for $\tilde \phi$ as given
in (2.22) and by simplifying, we can show that inequality (2.37)
is equivalent to $B_2 < B_1$, which is certainly true because $B_2
< 0$ and $B_1 > 1$.  Thus, we have proved the following theorem:

\th{2.3} {For $M < m \le w < x$, the optimal investment strategy
$\pi^\phi(w) > {\mu - r \over \sigma^2} {1 \over p - 1} \left( {c
\over r} - w \right).$}  $\square$

\medskip

Therefore, not only have we {\it not} reduced the leveraging by
considering the probability of ruin at death, we have made it
strictly worse for $M < m \le w < x$.  The following example shows
the extent of this worsening for some reasonable model parameters.

\ex{2.4} Assume the following parameter values for a model in
which consumption is expressed in real terms, that is, after
inflation:
\item{$\bullet$} $\lambda = 0.04$; the hazard rate is constant such that the expected future lifetime is 25 years.
\item{$\bullet$} $r = 0.02$; the riskless rate of return is 2\% over inflation.
\item{$\bullet$} $\mu = 0.06$; the risky asset's drift is 6\% over inflation.
\item{$\bullet$} $\sigma = 0.20$; the risky asset's volatility is 20\%.
\item{$\bullet$} $c = 1$; the individual consumes one unit of real wealth per year.

\item{$\bullet$} $x = 0$; the ruin level at death is 0.

\item{$\bullet$} $M = -200$; the lifetime ruin level is -200.

See Figure 1 for a graph of the optimal investment strategy
$\pi^\phi$ compared with the optimal investment strategy for the
problem of minimizing lifetime ruin, as in Section 2.1, namely $
{\mu - r \over \sigma^2} {1 \over p - 1} \left( {c \over r} - w
\right).$  Note the discontinuity of $\pi^\phi(w)$ when $w = x =
0$, as expected from the discontinuity of (2.10).  More
importantly, note that $\pi^\phi(w) > {\mu - r \over \sigma^2} {1
\over p - 1} \left( {c \over r} - w \right)$ for $w < 0$, as shown
in Theorem 2.3.

\medskip

\centerline{\bf Figure 1 about here}

\medskip

In the next section, we remove the leveraging entirely by
prohibiting borrowing of the riskless asset, as in Bayraktar and
Young (2007a).

\subsect{2.3. Probability of Lifetime Ruin with No Borrowing}

Bayraktar and Young (2007a) consider the problem of minimizing the
probability of lifetime ruin under the constraint that the
individual cannot borrow; however, they consider only the case for
which the ruin level $x = 0$.  In Section 3.3, we require the
corresponding minimum probability of ruin for an arbitrary level
of ruin $x < c/r$, so in this section, we consider that problem.

Let $\p(w, m; x)$ denote the minimum probability that the
individual's wealth reaches $x < c/r$ before she dies under the
constraint that she cannot borrow money to invest in the risky
asset.  We will consider two cases when wealth reaches zero: (1)
Welfare provides income at the rate of $c$ to provide for the
consumption needs of the individual.  In that case, zero is an
absorbing state for the wealth process, so if $W_t = 0$, then
$W_{t+s} = 0$ for $s \ge 0$.  (2) The individual is allowed to
borrow but only to cover her consumption, that is, she cannot
borrow to invest in the risky asset.  In this case if $W_t = 0$,
then $W_{t+s} = -c(e^{rs} - 1)/r$ for $s \ge 0$.

When $x = 0$, Bayraktar and Young (2007a) show that $\p$ is given
by

$$\p(w, m; 0) = \cases{1, &if $m \le 0$; \cr h_0(w), &if $0 < m
\le w \le w_l$; \cr \beta_0 \left( 1 - {r \over c} w \right)^p,
&if $0 < m$ and $w_l < w < c/r$; \cr 0, &if $0 < m$ and $w \ge
c/r$;} \eqno(2.38)$$

\noindent in which $w_l = {\xi \over 1 + \xi} {c \over r}$, where
$\xi = {\mu-r \over \sigma^2}{1 \over p - 1}$, and $h_0$ solves

$$\left\{ \eqalign{\lambda h_0 &= (\mu w - c) h_0' + {1 \over 2}
\sigma^2 w^2 h_0'', \cr h_0(0) &= 1, \quad {h_0(w_l) \over
h'_0(w_l)} = -{1 \over p} \left( {c \over r} - w_l \right).}
\right.  \eqno(2.39)$$

\noindent Once we have $h_0$, then we can compute $\beta_0 =
h_0(w_l)(1 - rw_l/c)^{-p}$.  The corresponding optimal investment
strategy for $0 < m \le w < c/r$ is given by

$$\pi^{nb}(w) = \cases{w, & if $0 < w \le w_l;$ \cr {\mu - r \over
\sigma^2} {1 \over p - 1} \left( {c \over r} - w \right), & if
$w_l < w < c/r$.} \eqno(2.40)$$

\noindent The value $w_l$ is called lending level because when
wealth is greater than $w_l$, the individual ``lends'' money to
the bank by buying the riskless asset.

We now extend Bayraktar and Young (2007a) to arbitrary $x < c/r$
for two reasons:  (1) The optimal investment strategy is
independent of the ruin level $x$, an interesting result in
itself.  (2) In Section 3.3, when we consider minimizing the
expected lifetime shortfall under a no borrowing constraint, then
we will use the $\p$ to represent the value function.  For the
sake of brevity, we simply state the minimum probability of ruin
$\p$ and the optimal investment strategy $\pi^{nb}$ because the
proof of these results are similar to those of Bayraktar and Young
(2007a).

\medskip

\noindent {\bf Case 1.} $w_l \le x < c/r$.  In this case, the
constraint will not bind, and we have

$$\p(w, m; x) = \cases{1, &if $m \le x$; \cr \left( {c - rw \over
c - rx} \right)^p, &if $x < m \le w < c/r$; \cr 0, &if $x < m$ and
$w \ge c/r$;} \eqno(2.41)$$

\noindent and the corresponding optimal investment strategy for $x
< m \le w < c/r$ is given by

$$\pi^{nb}(w) = {\mu - r \over \sigma^2} {1 \over p - 1} \left( {c
\over r} - w \right).  \eqno(2.42)$$

\noindent Note that the investment strategy in (2.42) agrees with
the one given in (2.40) on their common domain.

\medskip

\noindent {\bf Case 2.} $0 \le x < w_l$.  In this case, one can
parallel the argument in Bayraktar and Young to show that the
constraint does not bind for $w \in (w_l, c/r)$ and does bind for
$w \in (x, w_l)$.  Thus, we have

$$\p(w, m; x) = \cases{1, &if $m \le x$; \cr h_x(w), &if $x < m
\le w \le w_l$; \cr \beta_x \left( {c - rw \over c - rx}
\right)^p, &if $x < m$ and $w_l < w < c/r$; \cr 0, &if $x < m$ and
$w \ge c/r$;} \eqno(2.43)$$

\noindent in which $h_x$ solves (2.39) with the boundary condition
$h_0(0) = 1$ replaced by $h_x(x) = 1$.  Once we have $h_x$, then
we can compute $\beta_x = h_x(w_l)[(c - rw_l)/(c - rx)]^{-p}$.
The corresponding optimal investment strategy for $x < m \le w <
c/r$ is given by

$$\pi^{nb}(w) = \cases{w, & if $x < w \le w_l;$ \cr {\mu - r \over
\sigma^2} {1 \over p - 1} \left( {c \over r} - w \right), & if
$w_l < w < c/r$.} \eqno(2.44)$$

\noindent As in Case 1, the investment strategy in (2.44) agrees
with the one given in (2.40) on their common domain.

\medskip

\noindent {\bf Case 3.} $x < 0$.  In this case, we make one of two
assumptions as described in the paragraph preceding (2.38).  Under
the assumption that welfare will cover consumption if wealth
reaches zero, we have that zero is an absorbing state.  Also, if
the process starts with $w < 0$, the process will remain there
because of welfare.  The corresponding trivial observation is that
$\p(w, m; x) = 1$ if $m \le x$, and $\p(w, m; x) = 0$ if $m > x$.
There is no unique optimal investment strategy in this case.

Under the assumption that the individual is allowed to borrow to
cover consumption, we have that $\pi^{nb}(w) = 0$ for $w < 0$,
which immediately leads to the conclusion that

$$\p(w, m; x) = \cases{1, &if $m \le x$; \cr \left( {c - rw \over
c - rx} \right)^{\lambda/r}, &if $x < m \le w \le 0$; \cr h_x(w),
&if $x < m$ and $0 < w \le w_l$; \cr \beta_x \left( {c - rw \over
c - rx} \right)^p, &if $x < m$ and $w_l < w < c/r$; \cr 0, &if $x
< m$ and $w \ge c/r$;} \eqno(2.45)$$

\noindent in which $h_x$ solves (2.39) with the boundary condition
$h_0(0) = 1$ replaced by $h_x(0) = (c/(c - rx))^{\lambda/r}$.  The
corresponding optimal investment strategy for $x < m \le w < c/r$
is given by

$$\pi^{nb}(w) = \cases{0, & if $x < w \le 0$; \cr w, & if $0 < w
\le w_l;$ \cr {\mu - r \over \sigma^2} {1 \over p - 1} \left( {c
\over r} - w \right), & if $w_l < w < c/r$.} \eqno(2.46)$$

\medskip

In the next section, we parallel the work from this section with
an objective function that penalizes for the amount that an
individual's wealth falls below a given level, not just whether or
not wealth falls below this level.

\sect{3. Expected Lifetime Shortfall and Shortfall at Death}

Initially, we anticipated that if the individual is penalized by
the amount of loss, then she will take less chance with investing
in the risky asset than if she were penalized simply for her
wealth being low regardless of how low.  However, we were
surprised to learn that the investment strategy is the same when
minimizing expected lifetime shortfall as when minimizing the
probability of lifetime ruin, and we show in Section 3.1.

In Section 3.2, we minimize the expected shortfall at death and
show that the optimal investment strategy is larger than if
minimizing expected lifetime shortfall, in parallel to what we
learned in Section 2.2.  In other words, modifying the objective
function to penalize the individual for the amount of loss does
not ameliorate the leveraging effect, and considering wealth at
death only makes it worse.  Therefore, to reduce leveraging, in
Section 3.3, we eliminate it completely by prohibiting borrowing
in the problem in Section 3.1.

\subsect{3.1. Expected Lifetime Shortfall}

The individual's objective is to minimize the maximum shortfall
relative to $x < c/r$ during life. Then, the relevant value
function for this individual's objective is given by

$$V(w, m; x) = \inf_\pi \E[(x - M_{\tau_d})_+ \big| W_0 = w, M_0 =
m].  \eqno(3.1)$$

\noindent Here $(x - m)_+ = \max(x - m, 0)$ denotes the negative
part of the random variable $M_{\tau_d} - x$, and $\tau_d$ denotes
the random time of death of our individual, as in Section 2.  We
refer to $(x - M_{\tau_d})_+$ as the lifetime shortfall relative
to $x$.   The minimization in (3.1) is carried out over all
admissible investment strategies, as defined in Section 2.1.   In
this section and in Section 3.2, we apply no further constraints
on admissible investment strategies while in Section 3.3, we
require that $\pi_t \leq \max(0, W_t)$.

We first consider value functions of the form

$$V^f(w, m) = \inf_\pi \E[f(M_{\tau_d}) \big| W_0 = w, M_0 = m],
\eqno(3.2)$$

\noindent in which $f$ is non-increasing, bounded, and
continuously differentiable, and $f(m) = 0$ for $m \ge x$ for some
$x < c/r$. That is, $f$ is a non-decreasing function of shortfall
relative to $x$. We provide a verification theorem (see Appendix A
for its proof) and obtain (3.2) explicitly.  As a consequence of
this result, in Theorem 3.3, we find the optimal investment
strategy for the problem of minimizing the maximum expected
shortfall (3.1).  We also evaluate the value function in (3.1)
explicitly; see (3.19).

\th{3.1} {Let ${\bf D} = \{(w,m) \in {\bf R} \times {\bf R}: w
\geq m \}$. Suppose $h: {\bf D} \rightarrow {\bf R}$ is a bounded,
continuous function that satisfies the following conditions:
\item{$(i)$} $h(\cdot, m) \in C^2([m,c/r))$ for all $m \in {\bf R};$
\item{$(ii)$} $h(w, \cdot)$ is continuously differentiable$;$
\item{$(iii)$} $h_m(m, m)=0$ for all $m \in {\bf R};$
\item{$(iv)$} $h(w, m) = f(m)$ for $w \geq c/r;$
\item{$(v)$} $h_w(c/r, m)=0;$
\item{$(vi)$} $h$ solves the following $($HJB$)$ equation}
$$\lambda (h(w,m) - f(m)) = (rw - c)h_{w}(w, m) + \min_{\pi}
\left[(\mu-r)\pi h_{w}(w, m) + {1 \over 2}\sigma^2
\pi^{2}h_{ww}(w, m) \right]. \eqno(3.3)$$

\noindent{\it Then the value function in $(3.2)$ is given by}
$$V^f(w, m)=h(w, m), \quad -\infty < m \le w < \infty.
\eqno(3.4)$$

\medskip

We next use Theorem 3.1 to solve for $V^f$.  We hypothesize that
$V^f$ is convex (in its first variable), satisfies (i), and solves
the HJB equation (3.3) together with the boundary conditions (iv)
and (v) of Theorem 3.1.  Under these assumptions we obtain an
explicit expression for $V^f$ and later check that it satisfies
(ii) and (iii).

Since we assume that $V^f$ is convex, we can compute its concave
dual $\tilde V^f$ by the Fenchel-Legendre transform, as in the
proof of Theorem 2.1.  For ease of reference, we repeat the
analogs of (2.13) and (2.14) here.

$$\Vt^f (y, m) = \min_w [V^f(w, m)+ w y].  \eqno(3.5)$$

\noindent Note that we can retrieve the function $V^f$ from
$\Vt^f$ by the relationship

$$V^f(w,m) = \max_y [\Vt^f(y,m) - w y]. \eqno(3.6)$$

\noindent The minimizer of the right-hand side of (3.5) equals
$I(-y,m) = \Vt^f_y(y,m)$, where $I$ is the inverse function of
$V^f_w$ with respect to $w$.  Therefore, the maximizer of the
right-hand side of (3.6) is equal to $-V^f_{w}(w,m)$.

Substitute $w = I(-y,m)$ into (3.3) to get

$$\lambda \Vt^f (y, m)+ (r - \lambda)y \Vt^f_y (y, m) - \delta y^2
\Vt^f_{yy}(y, m) = cy + f(m), \eqno(3.7)$$

\noindent in which  $\delta$ is given in (2.7).  The general
solution of (3.7) is given by

$$\Vt^f(y, m) = D_1 y^{B_1} + D_2 y^{B_2} + {c \over r} y + f(m),
\eqno(3.8)$$

\noindent in which $B_1$ and $B_2$ are given in (2.18) and (2.19),
respectively.  $D_1$ and $D_2$ are functions of $m$ to be
determined.

Define $y_c = -V^f_w (c/r, m) = 0$; that is, $\Vt^f_y (0, m) =
c/r$. From the definition of $\Vt^f$ and from $V^f(c/r, m) =
f(m)$, we have

$$\Vt^f(0,m) = V^f(c/r, m) + {c \over r} y_c = f(m). \eqno(3.9)$$

\noindent Expression (3.9) implies that $D_2 = 0$ in (3.8).  Now,
we can calculate $V^f$ from

$$V^f(w, m) = \max_y \left[D_1 y^{B_1} + {c \over r} y + f(m) - wy
\right].  \eqno(3.10)$$

\noindent From the first-order condition, the maximizer of the
right-hand side of (3.10) is given by

$$y=\left(-{c/r - w \over D_1 B_1} \right)^{1/(B_1 - 1)}.
\eqno(3.11)$$

By substituting (3.11) into (3.10), we obtain

$$V^f(w,m) = k(m) \left({c \over r} - w \right)^p + f(m),
\eqno(3.12)$$

\noindent for some $k(m) > 0$ yet to be determined and for $p$
given in (2.6).

\noindent Now, we obtain an explicit expression for $\pv$ by using
the convexity of $V^f$ and the corresponding first-order condition
from (3.3), namely

$$\pv(w, m) = - {\mu - r \over \sigma^2} {h_{w}(w, m) \over
h_{ww}(w, m)}, \eqno(3.13)$$

\noindent and obtain that $\pv$ is {\it identical} to the
investment strategy in (2.8) for the problem of minimizing the
problem of lifetime ruin.

We next show how the minimum probability of lifetime ruin can be
used to compute $k(m)$ in (3.12).  Denote by $M^{*}$ the minimum
wealth process if the individual follows the optimal investment
strategy $\pv = \pi^\psi$. Then, by using the fact that the
optimal investment strategies are equal for $V^f$ and for $\psi$,
we can write

$$\P \left(M^*_{\tau_d} \le y \big| W_0 = w, M_0 = m \right) =
\cases {\left({c - rw  \over c - ry} \right)^p, &for $y < m <c/r$,
\cr 1, &for $m \le x < c/r$,} \eqno(3.14)$$

\noindent for $w < c/r$.  It follows that we can write (3.2) as

$$\eqalign{V^f(w, m) &= \E_{w,m}[f(M^*_{\tau_d}) \big] =
\int_{-\infty}^\infty f(y) d\P\left(M^*_{\tau_d} \le y \big| W_0 =
w, M_0 = m \right) \cr & = f(m) \left(1 - \psi(w, m; m) \right) +
\int_{-\infty}^{m}f(y) d\psi(w, m; y) \cr & = f(m) -
\int_{-\infty}^m f'(y) \left({c-r w \over c - ry} \right)^p dy,}
\eqno(3.15)$$

\noindent for $w < c/r$.  $\E_{w, m}$ denotes the conditional
expectation given $W_0 = w$ and $M_0 = m$.

Note that if $w \ge c/r$, then $\psi \equiv 0$, and $V^f(w, m) =
f(m)$, which is consistent with the first and second lines of
(3.15).  From (3.15) it follows that $V^f(w, \cdot)$ is
continuously differentiable because $f$ is continuously
differentiable and that $V^f_m(m, m)=0$. We summarize our findings
in the next theorem.

\th{3.2} {With no constraints on borrowing, $V^f$ defined by
$(3.2)$ is equal to $$V^f(w,m) = f(m) - \int_{-\infty}^m f'(y)
\left({c-r w \over c - ry} \right)^p dy, \eqno(3.16)$$ for $m \le
w < c/r$. The optimal investment strategy is given by $$\pv (w) =
{1 \over p-1} {\mu-r \over \sigma^2} \left({c \over r} - w
\right), \eqno(3.17)$$ for $w < c/r,$ which is equal to the
optimal investment strategy that attains the minimum probability
of lifetime ruin $\psi$.}  $\square$

\medskip

As a result of Theorem 3.2, the following theorem (whose proof is
in Appendix B) shows that $\pv$ in (3.17) is also an optimal
investment strategy to minimize the maximum lifetime shortfall.

\th{3.3} {Let $\pv$ be as in Theorem 3.2; then, the value function
of the shortfall problem $(3.1)$ satisfies $$V(w, m; x) = \E_{w,m}
\left[(x - M^{\pv}_{\tau_d})_+ \right]; \eqno(3.18)$$ that is,
$\pv$ minimizes the expected maximum lifetime shortfall.}

\medskip

As a corollary to Theorem 3.3, we can use the expression in (3.16)
to write (3.1) as

$$\eqalign{V(w, m; x) &= (x - m)_+ \left( 1 - \left( {c - rw \over
c - rm} \right)^p \right) + \int_{-\infty}^m (x - y)_+ {\partial
\over \partial x} \left( {c - rw \over c - ry} \right)^p dy \cr &=
(x - m)_+ + \left( {c - rw \over c - r(m \wedge x)} \right)^p {c -
r(m \wedge x) \over r(p - 1)}.} \eqno(3.19)$$

Besides the useful representations of $V^f$ and $V$ in (3.16) and
(3.19), respectively, the main fact to glean from this section is
that the optimal investment strategy for the problems in (3.1) and
(3.2) are identical and are identical to the one for minimizing
the probability of ruin, as given in (2.8).  Thus, introducing a
penalty for the amount that one's wealth drops below $x$ does {\it
not} reduce the leveraging at small wealth in the optimal
investment strategy.  On the other hand, the fact that the optimal
investment strategies are the same leads to the useful
representation of $V^f$ in (3.16).

In the next section, we modify the penalty in (3.2) by considering
the shortfall at death.

\subsect{3.2. Expected Shortfall at Death}

As we were motivated to consider the problem in Section 2.2, we
hope that by introducing a penalty for the shortfall at death, we
will reduce the leveraging in the optimal investment strategy.  In
Section 2.2, we were careful in modifying the definition of $\psi$
so that $\phi$ was convex.  We take the same care here.

Note that one can write $\phi$ from (2.9) as

$$\phi(w, m; x, M) = \inf_\pi {\bf  E} \left[ {\bf
1}_{\{W_{\tau_d} \le x \}} {\bf  1}_{\{M_{\tau_d} > M \}} + {\bf
1}_{\{M_{\tau_d} \le M \}} \Big| W_0 = w, M_0 = m \right].
\eqno(3.20)$$

\noindent By comparing (3.20) with $\psi$ in (2.4), specifically
with

$$\psi(w, m; x) = \inf_{\pi} {\bf  E} \left[{\bf  1}_{\{M_{\tau_d}
\le x \}} \vert W_0 = w, M_0 = m \right], \eqno(3.21)$$

\noindent we now modify the penalty function in (3.1) in a similar
manner and define $U$ as follows:

$$U(w, m; x, M) = \inf_\pi {\bf  E} \left[ \left( x - W_{\tau_d}
\right)_+ {\bf  1}_{\{M_{\tau_d} > M \}} + (x - M) {\bf
1}_{\{M_{\tau_d} \le M \}} \Big| W_0 = w, M_0 = m \right],
\eqno(3.22)$$

\noindent in which $M < x < c/r$.

For $m \le M$, $U$ is identically $x - M$, and for $m > M$ and $w
\ge c/r$, $U$ is identically 0.  For $M < m \le w < c/r$, $U$
solves the following HJB equation:

$$\left\{ \eqalign{& \lambda (U(w) - (x - w)_+) = (rw - c) U'(w) +
\min_\pi \left[ (\mu - r) \pi U'(w) + {1 \over 2} \sigma^2 \pi^2
U''(w) \right], \cr & U(M) = x - M, \quad U(c/r) = 0,} \right.
\eqno(3.23)$$

\noindent in which we drop the notational dependence of $U$ on $m$
because $U$ is independent of $m$ if $m > M$.

We have the following theorem whose proof we omit because it is
similar to that of Theorem 2.1.

\th{3.4} {For $m > M$ and $x \le w < c/r$, the function $U$ is
given by}

$$U(w, m; x, M) = \beta \left( 1 - {r \over c} w \right)^p,
\eqno(3.24)$$

\noindent{\it with $p$ as in $(2.6)$ and $\beta > 0$.  Thus, on
$[x, c/r)$, $U$ is a multiple of the probability of lifetime ruin
$\psi$, or equivalently, a multiple of $V - (x - m)_+$ from
$(3.19)$.  It follows that the optimal investment strategy $\pi^U$
equals $\pi^\psi$, as given in $(2.8)$.}  $\square$

\medskip

Next, we consider the solution of (3.23) for $M < m \le w < x$.
On this domain, $U$ solves the system

$$\left\{ \eqalign{& \lambda (U(w) - (x - w)) = (rw - c) U'(w) +
\min_\pi \left[ (\mu - r) \pi U'(w) + {1 \over 2} \sigma^2 \pi^2
U''(w) \right], \cr & U(M) = x - M, \quad U(x)/U'(x) = -(c/r -
x)/p,} \right. \eqno(3.25)$$

\noindent in which we assume that $U$ and $U'$ are continuous at
$w = x$.  As in Sections 2.2 and 3.1, we consider the
Fenchel-Legendre transform of $U$, denoted by $\tilde U$.  We show
how to solve for $\tilde U$ numerically, and from (2.14) or (3.6),
we can then determine $U$.  Even though we can only calculate $U$
numerically, we show that $\pi^U > \pi^\psi$ for $M < m \le w <
x$.

By substituting $w = I(-y, m)$ in (3.25), we obtain the following
equation:

$$\lambda \tilde U(y) + [(r - \lambda)y + \lambda] \tilde U'(y) -
\delta y^2 \tilde U''(y) = \lambda x + cy.  \eqno(3.26)$$

\noindent Write $\tilde U$ as the sum of a solution to the
homogeneous problem and of a particular solution:

$$\tilde U(y) = \tilde U_h(y) + {c \over r} y - \left( {c \over r}
- x \right).  \eqno(3.27)$$

If we define $y_x = -U'(x) = \beta p (r/c) (1 - rx/c)^{p-1}$, then
we have

$$\tilde U_h(y_x) = \left( {c \over r} - x \right) \left( 1 - {p -
1 \over p} y_x \right), \eqno(3.28)$$

\noindent and

$$\tilde U'_h(y_x) = x - {c \over r}.  \eqno(3.29)$$

\noindent  Similarly, if we define $y_M = -U'(M)$, then we have

$$\tilde U_h(y_M) = \left( {c \over r} - M \right) \left( 1 - y_M
\right), \eqno(3.30)$$

\noindent and

$$\tilde U'_h(y_M) = M - {c \over r}.  \eqno(3.31)$$

Therefore, the second-order linear differential equation for
$\tilde U_h$, together with the four equations (3.28)-(3.31) at
the two boundaries, are enough to determine $\tilde U_h$, $y_x$,
and $y_M$.

To end this section, we compare the optimal investment strategy
for this problem $\pi^U$ with the one for the problems in Sections
2.1 and 3.1, namely $\pi^\psi$ given in (2.8).  Recall that for $m
> M$ and $x \le w < c/r$, the two strategies are identical.  For
$M < m \le w < x$, $\pi^U(w) > \pi^\psi(w)$ if and only if

$$-y \tilde U''(y) > {1 \over p - 1} \left( {c \over r} - \tilde
U'(y) \right), \quad y_x < y < y_M. \eqno(3.32)$$

\noindent Compare this inequality with (2.37).  In terms of
$\tilde U_h$, we can rewrite (3.32) as follows:

$$y \tilde U''_h(y) < {1 \over p - 1} \tilde U'_h(y), \quad y_x <
y < y_M. \eqno(3.33)$$

Substitute for $\tilde U''_h$ in (3.33) from the homogeneous
equation derived from (3.27) and simplify to obtain the inequality

$${\tilde U_h(y) \over \tilde U'_h(y)} > {p - 1 \over p} y - 1,
\quad y_x < y < y_M.  \eqno(3.34)$$

\noindent Thus, if we can demonstrate that inequality (3.34)
holds, then $\pi^U(w) > \pi^\psi(w)$ for $M < m \le w < x$.
Inequality (3.34) does hold, and we demonstrate this in the proof
of the following theorem:

\th{3.5} {For $M < m \le w < x$, the optimal investment strategy
$\pi^U(w) > {\mu - r \over \sigma^2} {1 \over p - 1} \left( {c
\over r} - w \right).$}

\medskip

\pf Inequality (3.34) holds if and only if $g(y) > z(y)$ for $y_x
< y < y_M$, in which $g$ and $z$ are defined by $g(y) = \tilde
U_h(y)/ \tilde U'_h(y)$ and $z(y) = [(p - 1)/p] y - 1$.  For $0 <
y_x < y < y_M$, $g$ solves the following first-order non-linear
differential equation:

$$g'(y) = 1 - {\lambda \over \delta y^2} g^2(y) - {(r - \lambda) y
+ \lambda \over \delta y^2} g(y), \eqno(3.35)$$

\noindent with boundary condition at $y = y_x$ determinable from
(3.28) and (3.29).  Specifically, $g(y_x) = z(y_x) =  [(p - 1)/p]
y_x - 1$.

Suppose $g(y) = z(y)$ for some $y_x < y < y_M$, then from (3.35),
we deduce that $g'(y) > z'(y) = (p - 1)/p$ if and only if $rp >
\lambda$, independent of $y$.  Now, $rp > \lambda$ for $r < \mu$;
therefore, if $g$ and $z$ intersect at some point in $[y_x, y_M)$,
then the slope of $g$ is larger than the slope of $z$ at that
point.  In particular, we know that $g(y_x) = z(y_x)$; therefore,
$g'(y_x) \ge z'(y_x)$.  In order for $g$ to intersect $z$ at $y >
y_x$, we must have $g'(y) \le z'(y)$ at the smallest such value of
$y$, a contradiction to $g'(y) > z'(y)$ whenever $g(y) = z(y)$.
Therefore, no such intersection point $y > y_x$ exists, and $g$ is
strictly larger than $z$ on $(y_x, y_M)$.  $\square$

\medskip

From Theorem 3.5, we learn that as in Section 2.2, not only did we
not reduce the leveraging exhibited in the investment strategy
given in (2.8), we actually increased the leverage by considering
a penalty that depends on the shortfall at death.  Therefore, as
in Section 2, our only real hope of reducing leverage is to
prevent it by modifying the admissible investment strategies.  We
do that in the next section, in parallel to Section 2.3.

\subsect{3.3. No Borrowing}

In this section, we consider the problem of minimizing the
expected maximum lifetime shortfall when the individual is not
allowed to borrow money to invest in the stock market. We will
observe that the optimal investment strategy that minimizes the
probability of ruin (independent of the ruin level $x$) and the
optimal investment strategy that minimizes the expected maximum
shortfall are the same, subject to the constraint that $\pi_t \le
\max(0, W_t)$.  As in Section 2.3, we consider two assumptions
concerning negative wealth: (1) Welfare provides for consumption,
and (2) the individual can borrow only to cover consumption.

As in Section 3.1, we first consider a continuously
differentiable, bounded, non-increasing function $f$ of the
minimum wealth such that $f(m) = 0$ for $m > x$ for some $x <
c/r$. Without abusing the notation too much, we denote the
corresponding value function by $V^f$, as in (3.2).  We provide a
verification lemma whose proof we omit because of its similarity
to the proof of Theorem 3.1.  Recall that ${\bf D} = \{(w,m) \in
{\bf R} \times {\bf R}: w \geq m \}$.

\th{3.6} {Suppose $a:{\bf D} \rightarrow {\bf R}$ is a bounded
continuous function that satisfies the following conditions:
\item{$(i)$} $a(\cdot, m) \in C^{2}([m, c/r))$ for all $m \in {\bf R}$;
\item{$(ii)$} $a(w, \cdot)$ is continuously differentiable$;$
\item{$(iii)$} $a_m(m, m)=0$ for all $m \in {\bf  R};$
\item{$(iv)$} $a(w, m) = f(m)$ for $w \geq c/r;$
\item{$(v)$} $a_w(c/r, m)=0;$
\item{$(vi)$} $a$ solves the following HJB equation for $0 < w < c/r$}
$$\lambda(a(w,m) - f(m)) = (rw - c) a_w(w, m) + \min_{\pi \le w}
\left[ (\mu - r)\pi a_w(w, m) + {1 \over 2} \sigma^2 \pi^2
a_{ww}(w, m) \right]. \eqno(3.36)$$
\item{({\it vii})} {\it If welfare exists, then $a(w, m) = f(m)$ for $w \le 0$.  If the individual is allowed to borrow to cover consumption when wealth is negative, then for $w \le 0$, the function $a$ solves $(3.36)$ with $\pi = 0$.}

\noindent {\it Then, the minimum expected maximum shortfall is
given by} $$V^f(w, m) = a(w, m), \quad -\infty < m \leq w <
\infty. \eqno(3.37)$$  $\square$

\medskip

Note that if welfare exists, then $V^f(w, m) = f(m)$ for $w \le
0$.  Alternatively, if the individual is allowed to borrow to
cover consumption, then $\pi(w) = 0$ for $w \le 0$, as in the case
for minimizing the probability of ruin in Section 2.3.  Thus, for
$w \le 0$,

$$\P(M^*_{\tau_d} \le y | W_0 = w, M_0 = m) = \cases{1, &if $m \le
y$; \cr \p(w, m; y) = \left( {c - rw \over c - rx}
\right)^{\lambda/r}, & if $m > y$.} \eqno(3.38)$$

\noindent It follows from (3.38) that we can write $V^f$ for $w
\le 0$ as follows:

$$\eqalign{V^f(w, m) &= \E[f(M^*_{\tau_d} | W_0 = w, M_0 = m] =
\int_{-\infty}^\infty d\P(M^*_{\tau_d} \le y | W_0 = w, M_0 = m)
\cr & = f(m) - \int_{-\infty}^m f'(x) \left( {c - rw \over c - ry}
\right)^{\lambda/r} dy.} \eqno(3.39)$$

Note that (3.39) gives us a boundary value for $V^f$ at $w = 0$
when borrowing is allowed for covering consumption.  One can apply
the arguments in Bayraktar and Young (2007a) to the more general
problem outlined in Theorem 3.6 and, thereby, we prove the
following theorem:

\th{3.7} {If welfare exists, then $V^f(w, m) = f(m)$ for $m \le w
\le 0$.  Alternatively, if the individual is allowed to borrow to
cover consumption when wealth is negative, then $V^f$ is given in
$(3.39)$ for $m \le w \le 0$.}

{\it For $0 < w \le w_l,$ $V^f$ solves}

$$\lambda(\V(w, m) - f(m)) = (\mu w - c)\V_w(w, m) + {1\over
2}\sigma^2 w^2 \V_{ww}(w, m), \eqno(3.40)$$ {\it with boundary
condition at $w = 0$ given by $V^f(0, m) = f(m)$ if there is
welfare or by $V^f(0, m)$ in $(3.39)$ if borrowing is allowed to
cover consumption and with boundary condition at $w = w_l$}
$${\V(w_l, m) - f(m) \over \V_w(w_l, m)} = -{1 \over p} \left({c
\over r} - w_l \right). \eqno(3.41)$$

{\it For $w_l < w < c/r,$ $V^f(w, m) = k(m) \left({c \over r} - w
\right)^p + f(m),$ in which $k(m) = \V(w_l, m)(1 - r w_l/c)^{-p}$.
Moreover, $V^f(w, m) = f(m)$ for $w \ge c/r$.}

{\it If borrowing is allowed to cover consumption when wealth is
negative, then for all $x < c/r,$ the optimal investment strategy
for the probability of ruin with ruin level $x$ and the optimal
investment strategy for $V^f$ coincide on $[x, c/r)$. If welfare
exists, then the optimal investment strategy for $\V$ when wealth
is positive is given by $(2.40)$.} $\square$

\medskip

We have a corollary that follows from the observation in Theorem
3.7 that the optimal investment strategy for the problem of
minimizing lifetime shortfall is the same as the strategy for
minimizing the probability of lifetime ruin.  Recall that we
assume that $f(m) = 0$ for $m \ge x$, in which $x$ is some number
less than $c/r$.

\cor{3.8} {If borrowing is allowed to cover consumption when
wealth is negative, then for $w < c/r,$ $$\P(\Ma \le y |W_0 = w,
M_0 = m) = \cases{\psi^{nb}(w, m; y), &for $y < m < c/r,$ \cr 1,
&for $m \le y,$} \eqno(3.42)$$}

\noindent{\it in which $\psi^{nb}$ is given in Cases 1-3 in
Section 2.3 and $w < c/r$. Therefore, if borrowing is allowed to
cover consumption,}

$$\V(w, m) = f(m) - \int_{-\infty}^m f'(y) \psi^{nb}(w, m; y) dy.
\eqno(3.43)$$

{\it If welfare exists to cover consumption when wealth is
negative, then}

$$\P(\Ma \le y |W_0 = w, M_0 = m) = \cases{0, &for $y < m \le w
\le0$, \cr 0, &for $y < \min(0, m)$ and $w > 0$, \cr \psi^{nb}(w,
m; y), &for $0 \le y < m < c/r,$ \cr 1, &for $m \le y,$}
\eqno(3.44)$$

\noindent{\it in which $\psi^{nb}$ is given in Cases 1 and 2 in
Section 2.3 and $w < c/r$.  Therefore, if welfare exists,}

$$\V(w, m) = \cases{f(m), &if $m \le 0,$ \cr f(m) - \int_0^m f'(y)
\psi^{nb}(w, m; y) dy, &if $m > 0$.} \eqno(3.45)$$

\medskip

\pf We prove (3.45) because (3.43) follows from (3.42) as (3.15)
follows from (3.14). If $m \le 0$, then $\P(\Ma \le y |W_0 = w,
M_0 = m) = 0$ if $y < m$ and equals 1 if $y \ge m$.  Thus, if $m
\le 0$, $V^f(w, m) = \int_{-\infty}^\infty f(y)
d\P\left(M^*_{\tau_d} \le y \big| W_0 = w, M_0 = m \right) =
f(m)(1 - 0) = f(m)$.

If $m > 0$, then

$$\eqalign{V^f(w, m) &= \int_{-\infty}^\infty f(y)
d\P\left(M^*_{\tau_d} \le y \big| W_0 = w, M_0 = m \right) \cr &=
f(0) \psi^{nb}w, m; 0) + \int_0^m f(y) d\psi^{nb}(w, m; y) +
f(m)(1 - \psi^{nb}(w, m; m)) \cr & = f(m) - \int_0^m f'(y)
\psi^{nb}(w, m; y) dy,} \eqno(3.46)$$

\noindent in which the last equality follows from integration by
parts.  $\square$

\medskip

As in Theorem 3.3, we can show that the optimal investment
strategy when $f(m) = (x - m)_+$ is identical to that in the
problem of minimizing the probability of lifetime ruin.  Denote
the corresponding value function by $V^{nb}$.  Therefore, we have
the following theorem:

\th{3.9} {If borrowing is allowed to cover consumption when wealth
is negative, then for $f(m) = (x - m)_+,$ with $x < c/r,$ the
value function of the constrained shortfall problem is given by}

$$V^{nb}(w, m; x) = (x - m)_+ + \int_{-\infty}^{m \wedge x}
\psi^{nb}(w, m; y)dy. \eqno(3.47)$$

\noindent {\it If welfare exists to cover consumption when wealth
is negative, then}

$$V^{nb}(w, m; x) = \cases{(x - m)_+, &if $m \le 0,$ \cr (x - m)_+
+ \int_0^{m \wedge x} \psi^{nb}(w, m; y) dy, &if $m > 0$.}
\eqno(3.48)$$ $\square$

\sect{4. Summary and Future Research}

In this paper, we examined the problems of minimizing (1) lifetime
ruin probability, (2) ruin probability at death, (3) expected
lifetime shortfall, and (4) expected shortfall at death. We showed
that the leveraging effect, the fact that an individual borrows
excessively at low wealth levels observed by Young (2004), is not
reduced by the alternative penalty criteria. In fact, we learned
that the optimal investment strategies of (1) and (3) are the
same, and that the optimal amount of wealth traded in the risky
asset for the cases of (2) and (4) exacerbate the leveraging
observed in (1). We also introduced a no-borrowing to constraint
to (1) and (3) to eliminate leveraging completely.

In work not shown here (for the sake of brevity), we also
considered the penalty functions (3) in a setting where we allowed
the individual to borrow but only at a rate greater than the one
earned by the riskless asset and did not impose a constraint on
the individual's investment strategy.  We found out that the
leveraging effect is exacerbated in some cases, and that the
optimal investment strategy is the same as in the corresponding
problem for minimizing the probability of lifetime ruin, which is
given in Bayraktar and Young (2007a).  We also observed that
changing the consumption function to a piece-wise linear function,
$c(w) = c + p(w-d)_+$, does not alter the leveraging effect.
Therefore, the leveraging effect is not a side effect of choosing
a constant consumption.

The leverage effect might decrease if one were to introduce
negative unbounded jumps, as found in Liu, Longstaff, and Pan
(2003) in the context of maximizing expected utility of terminal
wealth.  In future research, we plan to consider the problem of
minimizing the probability of lifetime ruin with stocks that are
subject to such negative jumps.

Also, in the future, instead of prohibiting borrowing entirely, we
will consider constraints of the from $\pi_t \leq k (W_t+ L)$.
Here $k, L > 0$.  Note that the constraint we considered in this
paper corresponds to $k = 1$ and $L = 0$. In general, for $k=1$,
the quantity $L$ might be the investor's credit limit.  The case
for which $L = 0$ and $0 < k < 1$ could represent the situation
for which the individual is allowed to put only a fraction of his
wealth into the risky asset.

\sect{Appendix A. Proof of Theorem 3.1}

Assume that $h$ satisfies the conditions specified in the
statement of Theorem 3.1. Let $N$ denote a Poisson process with
rate $\lambda$ that is independent of the standard Brownian motion
$B$ driving the wealth process. The occurrence of a jump in the
Poisson process represents the death of the individual investor.

Let $\pi: {\bf D} \rightarrow {\bf R}$ be a function, and let
$W^\pi$ and $M^\pi$ denote the wealth and the minimum wealth
respectively, when the individual investor uses the investment
policy $\pi_t = \pi(W_t, M_t)$.  Assume that this investment
policy is admissible.  For a given $m$, define ${\bf D}^m =
\{(w,m): w \geq m\}$, and define ${\bf \bar{D}}^m = {\bf D}^m \cup
\{\infty\}$ to be the one-point compactification of ${\bf D}^m$.
The point $\infty$ is  the ``coffin state." The wealth process is
killed (and sent to the coffin) as soon as the Poisson process
jumps (that is, when the individual dies), and we assign
$W^\pi_{\tau_d} = \infty$. All functions $g$ on ${\bf D}^{m}$ are
extended to ${\bf \bar{D}}^{m}$ by $g(\infty,m)= f(m)$. Observe
that $h(c/r,m)= h(W^{\pi}_{\tau_d},m)= f(m)$ for all $m \leq c/r
$.  Define the stopping time $\tau = \tau_d \wedge \tau_{c/r}$,
where $\tau_{c/r} = \inf\{t > 0: W_t = c/r \}$, with the
convention that $\inf \emptyset = \infty$.

If $w \geq c/r$, then the individual can invest her wealth in the
riskless asset and guarantee to finance the cost of her
consumption, which is almost surely less than $\int_{0}^{\infty} c
e^{-rt} dt = c/r$. Therefore, with this strategy, her wealth will
almost surely be at least $c/r$ at the time of her death. This
implies that $M_{\tau_d} = m$, and therefore $V^f(w, m)= f(m)$ for
$w \geq c/r$.  From this it follows that

$$V^f(w, m)=\inf_\pi \E[f(M_{\tau}) \big| W_0 = w, M_0 = m].
\eqno({\rm A.1})$$

Define

$$\tau_n = \inf \{t \ge 0: W_t \ge n \,\, \hbox{or}\,\,W_t \le
-n\,\, \hbox{or} \int_0^t \pi^2_s ds = n \}.  \eqno({\rm A.2})$$

\noindent By applying It\^{o}'s formula to $h$, we have

$$\eqalign{&h(W^\pi_{t \wedge \tau \wedge \tau_n}, M^\pi_{t \wedge
\tau \wedge \tau_n}) = h(w, m) \cr & \quad +\int_0^{t \wedge \tau
\wedge \tau_n }\left((rW^\pi_s + (\mu-r) \pi_s - c) h_w(W^\pi_s,
M^\pi_s) + {1 \over 2} \sigma^2 \pi_s^2 h_{ww}(W^\pi_s, M^\pi_s)
\right)ds \cr & \quad + \lambda \int_0^{t \wedge \tau \wedge
\tau_n} (f(M^\pi_s) - h(W^\pi_s, M^\pi_s))ds + \int_0^{t \wedge
\tau \wedge \tau_n} h_w(W^\pi_s, M^\pi_s) \sigma \pi_s dB_s \cr &
\quad + \int_0^{t \wedge \tau \wedge \tau_n} (f(M^\pi_s) -
h(W^\pi_s, M^\pi_s)) d(N_s - \lambda s) + \int_0^{t \wedge \tau
\wedge \tau_n} h_m (W^\pi_s, M^\pi_s) dM^\pi_s.} \eqno({\rm
A.3})$$

The last integral in (A.3) is equal to zero almost surely because
$dM_t$ is non-zero only when $M_t = W_t$, and $h_m(m,m) = 0$;
$h_m$ denotes the left derivative of $h$ with respect to $m$.
Here we also used the fact that $M$ is non-decreasing, therefore
the first variation process associated with it is finite almost
surely, to conclude that the cross variation of $M$ and $W$ is
zero almost surely.  It follows from the definition of $\tau_n$
that

$$\E_{w, m} \left[\int_0^{t \wedge \tau \wedge \tau_n} h_w
(W^\pi_s, M^\pi_s) \sigma \pi_s dB_s \right] = 0. \eqno({\rm
A.4})$$

\noindent $\E_{w, m}$ denotes the conditional expectation given
$W_0 = w$ and $M_0 = m$.  Moreover, the expectation of the fourth
integral is zero since $f$ and $h$ are bounded; see, for example,
Br\'emaud (1981).

Now, we have

$$\eqalign{& \E_{w,m} [h(W^\pi_{t \wedge \tau \wedge \tau_n},
M^\pi_{t \wedge \tau\wedge\tau_n})] = h(w,m) \cr & \quad + \E_{w,
m} \left[ \int_0^{t \wedge \tau \wedge \tau_n} \left((rW^\pi_s +
(\mu-r) \pi_s - c) h_w (W^\pi_s, M^\pi_s) + {1 \over 2} \sigma^2
\pi_s^2 h_{ww} (W^\pi_s, M^\pi_s) \right) ds \right] \cr & \quad +
\E_{w, m} \left[ \int_0^{t \wedge \tau \wedge \tau_n} \lambda
(f(M^\pi_s) - h(W^\pi_s, M^\pi_s)) ds \right] \cr & \geq h(w,m),}
\eqno({\rm A.5})$$

\noindent where the inequality follows from assumption (vi) of the
theorem.  Because $h$ is bounded by assumption, it follows from
the Dominated Convergence Theorem that

$$\E_{w,m}[h(W^\pi_{t \wedge \tau}, M^\pi_{t \wedge \tau})] \geq
h(w,m).  \eqno({\rm A.6})$$

\noindent Equation (A.6) shows that $(h(W^\pi_{t \wedge \tau},
M^\pi_{t \wedge \tau}))$ is a sub-martingale for any admissible
strategy $\pi$.

Since $h(c/r, m) = h(W^\pi_{\tau_d}, m) = f(m)$ for all $m \leq
c/r$, we have

$$h(W^\pi_{\tau}, M^\pi_{\tau}) = f(M_\tau). \eqno({\rm A.7})$$

\noindent If $\tau_d < \tau_{c/r}$, then obviously $M_\tau =
M_{\tau_d}$. If $\tau_d \ge \tau_{c/r}$, then $M_\tau =
M_{\tau_{c/r}}$.  By taking the expectation of both sides of
(A.7), we obtain

$$\E_{w,m} \left[h(W^\pi_{\tau}, M^\pi_{\tau}) \right] = \E_{w,m}
[f(M^\pi_\tau)] \geq h(w,m). \eqno({\rm A.8})$$

\noindent The inequality in (A.8) follows from an application of
the Optional Sampling Theorem because $(h(W^\pi_{t \wedge \tau},
M^\pi_{t \wedge \tau}))$ is a sub-martingale and $\sup_{t \geq 0}
\E_{w, m} [h(W^\pi_{t \wedge\tau}, M^\pi_{t \wedge \tau})] <
\infty$ because $h$ is bounded; see Theorem 3.15, page 17 and
Theorem 3.22, page 19 of Karatzas and Shreve (1991). Together with
(A.1) this implies that

$$V^f(w,m)=\inf_\pi \E[f(M^{\pi}_{\tau})] \geq h(w,m).  \eqno({\rm
A.9})$$

If the individual investor follows a strategy $\pv$ that minimizes
the right-hand side of (3.3), then (A.5) is satisfied with
equality, and an application of Dominated Convergence Theorem
yields

$$\E_{w,m} [h(W^{\pv}_{t \wedge \tau},M^{\pv}_{t \wedge
\tau})]=h(w,m), \eqno({\rm A.10})$$

\noindent which implies that $(h(W^{\pv}_{t \wedge \tau},
M^{\pv}_{t \wedge \tau}))$ is a martingale.  By following the same
line of argument as above, we obtain

$$V^f(w,m) = h(w,m), \eqno({\rm A.11})$$

\noindent which demonstrates that (3.4) holds and $\pv$ is an
optimal investment strategy.  $\square$

\sect{Appendix B. Proof of Theorem 3.3}

Let us introduce a sequence of increasing functions $g_n$ that
converges to $g(y) = (x - y)_+$, such that each element in this
sequence is bounded and is a difference of two convex functions:

$$g_n(y) = \cases{n, &if $y \le x - n$, \cr x - y, &if $x - n < y
\le x$, \cr 0, &if $y > x$.}  \eqno({\rm B.1})$$

\noindent  Note that $g_{n}(y) = (x - y)_+ - (x - n - y)_+$, and
$g_n(y) \leq n$. Observe that the sequence (B.1) indeed converges
monotonically to $g(y) = (x - y)_+$ as $n \rightarrow \infty$.
Note that $g_n$ is not continuously differentiable. We will
approximate $g_n$ with an increasing sequence $\tilde g_{n, k}$
given by

$$\tilde g_{n, k}(y) = \cases{n, &if $y  \le x - n - {1 \over k}$,
\cr n - {k \over 2} \left(y - x + n + {1\over k} \right)^2 , &if
$x - n - {1 \over k} < y \le x - n$, \cr - {1 \over 2k} - (y - x),
&if $x - n < y \le x - {1 \over k}$, \cr {1 \over 2k} - {1-k (y -
x) \over 2} \left( y - x + {1 \over k} \right), &if $x - {1 \over
k} < y \le x$, \cr 0, &if $y > x$.}  \eqno({\rm B.2})$$

\noindent The sequence $(\tilde g_{n, k})_{k \in {\bf N}}$
increases to $g_n$ for all $n$. Each element of the sequence of
$(\tilde g_{n, k})_{k \in {\bf N}}$ is continuously
differentiable.

We have

$$\eqalign{\E_{w,m} \left[(x - M^{\pv}_{\tau_d})_+ \right] &=
\E_{w,m} \left[ \lim_n g_n (M^{\pv}_{\tau_d}) \right] = \lim_n
\E_{w, m} \left[ g_n (M^{\pv}_{\tau_d}) \right] \cr & = \lim_n
\lim_k \E_{w,m} \left[ \tilde g_{n,k} (M^{\pv}_{\tau_d}) \right] =
\lim_n \lim_k \inf_\pi \E_{w, m} \left[\tilde
g_{n,k}(M^\pi_{\tau_d}) \right] \cr & \le \lim_n \E_{w,m}
\left[g_n(M^\pi_{\tau_d}) \right] \le \E_{w,m} \left[(x -
M^\pi_{\tau_d})_+ \right].}  \eqno({\rm B.3})$$

\noindent The second and the third equalities in (B.3) follow from
the Monotone Convergence Theorem, and the fourth equality follows
from Theorem 3.2. The inequalities follow from the fact that
$\tilde g_{n,k}(y) \le g_n(y) \le (x - y)_+$. By taking an infimum
over the admissible strategies, we obtain

$$\E_{w,m}\left[(x - M^{\pv}_{\tau_{d}})_+ \right] \le \inf_\pi
\E_{w,m} \left[(x - M^\pi_{\tau_d})_+ \right].  \eqno({\rm B.4})$$

\noindent This proves that $\pv$ in (3.17) is an optimal
investment strategy.  $\square$

\medskip

\centerline{\bf Acknowledgements} \medskip This research of the
first author is supported in part by the National Science
Foundation under grant DMS-0604491. We thank an anonymous referee,
Moshe A. Milevsky, Kristen S. Moore, and S. David Promislow for
their helpful comments.  We especially thank Thomas Salisbury for
suggesting this problem.

\sect{References}

\noindent \hangindent 20 pt Bayraktar, E. and V. R. Young (2007a),
Minimizing the probability of lifetime ruin under borrowing
constraints, to appear in {\it Insurance: Mathematics and
Economics.}

\noindent \hangindent 20 pt Bayraktar, E. and V. R. Young (2007b),
Correspondence between Lifetime Minimum Wealth and Utility of
Consumption, to appear in {\it Finance and Stochastics.}

\smallskip \noindent \hangindent 20 pt Br\'emaud, P. (1981), {\it Point Processes and Queues}, Springer Series in Statistics, Springer-Verlag, New York.

\smallskip \noindent \hangindent 20 pt Browne, S. (1995), Optimal investment policies for a firm with a random risk process: Exponential utility and minimizing the probability of ruin, {\it Mathematics of Operations Research}, 20 (4): 937-958.

\smallskip \noindent \hangindent 20 pt Browne, S. (1997), Survival and growth with liabilities: Optimal portfolio strategies in continuous time,  {\it Mathematics of Operations Research}, 22 (2): 468-493.

\smallskip \noindent \hangindent 20 pt Browne, S. (1999a), Beating a moving target: Optimal portfolio strategies for outperforming a stochastic benchmark, {\it Finance and Stochastics}, 3: 275-294.

\smallskip \noindent \hangindent 20 pt Browne, S. (1999b), Reaching goals by a deadline: Digital options and continuous time active portfolio management, {\it Advances in Applied Probability}, 31: 551-577.

\smallskip \noindent \hangindent 20 pt Browne, S. (1999c), The risks and rewards of minimizing shortfall probability, {\it Journal of Portfolio Management}, 25 (4): 76-85.

\smallskip \noindent \hangindent 20 pt Karatzas, I. and S. Shreve (1991), {\it Brownian Motion and Stochastic Calculus}, second edition, Springer-Verlag, New York.

\smallskip \noindent \hangindent 20 pt Karatzas, I. and S. Shreve (1998), {\it Methods of Mathematical Finance}, Springer-Verlag, New York.

\smallskip \noindent \hangindent 20 pt Liu, J., F. Longstaff, and J. Pan (2003), Dynamic asset allocation with event risk, {\it Journal of Finance}, 58 (1): 231-259.

\smallskip \noindent \hangindent 20 pt Merton, R. C. (1992), {\it Continuous-Time Finance}, revised edition, Blackwell Publishers, Cambridge, Massachusetts.

\noindent \hangindent 20 pt \O ksendal, B. (1998), {\it Stochastic
Differential Equations: An Introduction with Applications}, fifth
edition, Springer-Verlag, Berlin.

\smallskip \noindent \hangindent 20 pt Parikh, A. N. (2003), The evolving U.S. retirement system, {\it The Actuary}, March: 2-6.

\smallskip \noindent \hangindent 20 pt Pestien, V. C. and W. D. Sudderth (1985), Continuous-time red and black: How to control a diffusion to a goal, {\it Mathematics of Operations Research}, 10 (4): 599-611.

\smallskip \noindent \hangindent 20 pt Promislow, D. and V. R. Young (2005), Minimizing the probability of ruin when claims follow Brownian motion with drift, {\it North American Actuarial Journal}, 9 (3): 109-128.

\smallskip \noindent \hangindent 20 pt Roy, A. D. (1952), Safety first and the holding of the assets, {\it Econometrica}, 20: 431-439.

\smallskip \noindent \hangindent 20 pt Taksar, M. I. and C. Markussen (2003), Optimal dynamic reinsurance policies for large insurance portfolios, {\it Finance and Stochastics}, 7: 97-121.

\smallskip \noindent \hangindent 20 pt Schmidli, H. (2001), Optimal proportional reinsurance policies in a dynamic setting,  {\it Scandinavian Actuarial Journal}, 2001 (1): 55-68.

\smallskip \noindent \hangindent 20 pt Young, V. R. (2004), Optimal investment strategy to minimize the probability of lifetime ruin, {\it North American Actuarial Journal}, 8 (4): 105-126.

\smallskip \noindent \hangindent 20 pt Zariphopoulou, T. (1999), Transaction costs in portfolio management and derivative pricing, {\it Introduction to Mathematical Finance}, D. C. Heath and G. Swindle (editors), American Mathematical Society, Providence, RI. {\it Proceedings of Symposia in Applied Mathematics}, 57: 101-163.

\smallskip \noindent \hangindent 20 pt Zariphopoulou, T. (2001), Stochastic control methods in asset pricing, {\it Handbook of \break Stochastic Analysis and Applications}, D. Kannan and V. Lakshmikantham (editors), Marcel Dekker, New York.

\vfill \eject

\centerline{\vbox{\hsize=4.5in\psfig{figure=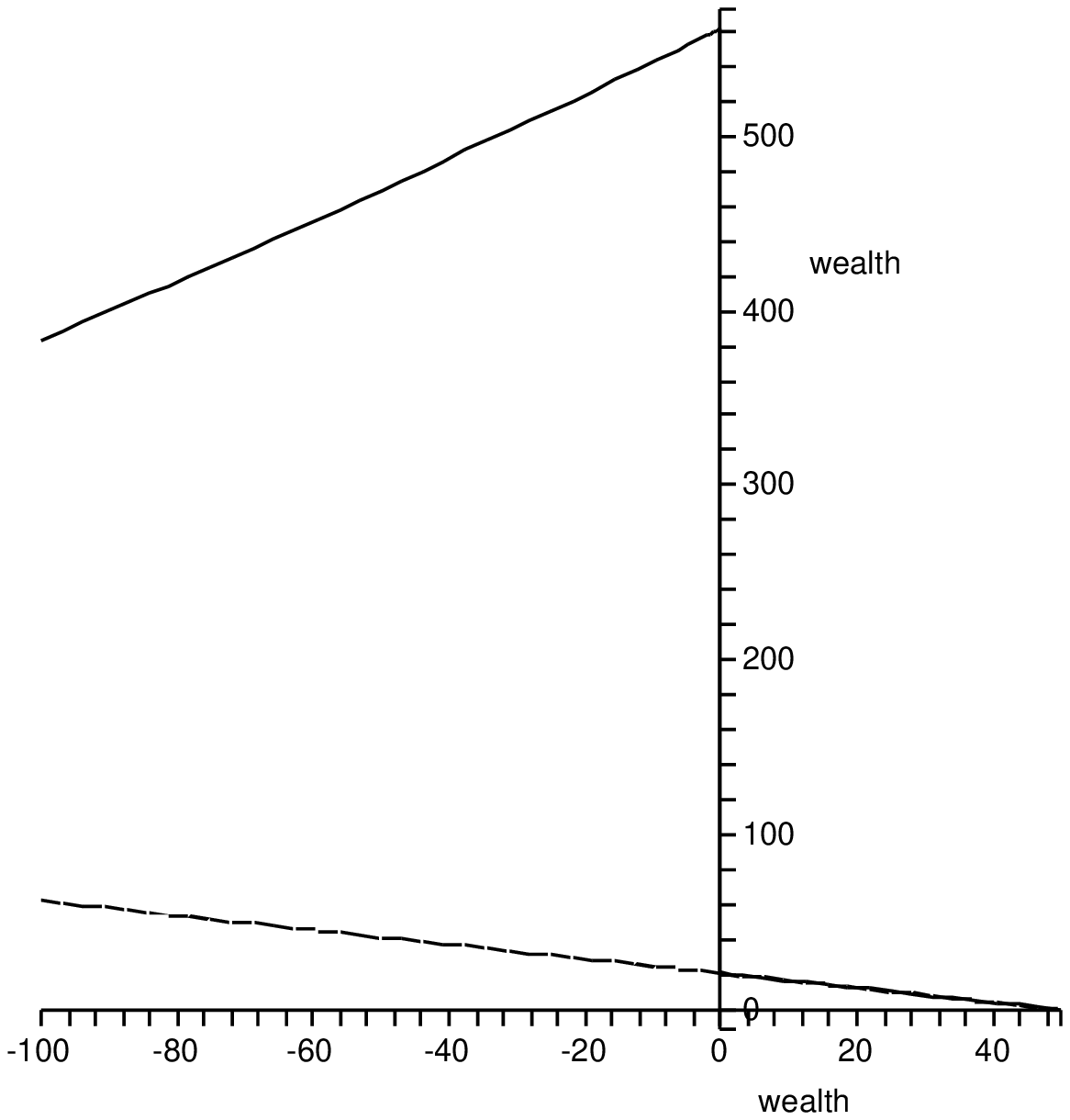,width=4.5in}
   \hfill \break
   Figure 1:  Graph of $\pi^\phi$ and $\pi^\psi$.  The solid line corresponds to $\pi^\phi$ and the dashed line to $\pi^\psi$.}}

\bye